\documentclass[12pt,notitlepage,a4paper,reqno,ps]{amsart}
\usepackage{eurosym}
\usepackage{amssymb}
\usepackage{amstext}
\usepackage{color}

\newcommand{\loc}{\textnormal{loc}}
\newcommand{\medint}{-\kern  -,375cm\int}

\newenvironment{michelarev}{\color{red}}{\color{black}}
\newcommand{\bmicr}{\begin{michelarev}}
\newcommand{\emicr}{\end{michelarev}}
\theoremstyle{plain}
\newtheorem{theorem}{Theorem}[section]

\theoremstyle{definition}

\theoremstyle{remark}
\newtheorem{remark}[theorem]{Remark}

\numberwithin{equation}{section} \makeatletter

\renewcommand{\p@enumi}{\thesection.}
\makeatother  \allowdisplaybreaks
\email{Corresponding author: mcaselli@student.ethz.ch (Michele Caselli)}

\begin{document}
\title[Lipschitz continuity results for obstacle problems]{Lipschitz continuity results \\for a class of obstacle problems }
\author[C. Benassi, M. Caselli]{Carlo Benassi, Michele Caselli}
\address{Dipartimento di Scienze Fisiche, Informatiche e Matematiche,
via Campi 213/b, 41125 Mo-dena, Italy}

\thanks{
The authors are partially supported GNAMPA (Gruppo Nazionale per l'Analisi Matematica, la Probabilit\`a e le loro Applicazioni)
 of INdAM (Istituto Nazionale di Alta Matematica) and by the University of Modena and Reggio Emilia through the project FAR2017 ``Equazioni differenziali: problemi evolutivi, variazionali ed applicazioni'' (coord. Prof. S.~Gatti). Part of the research has been performed while the second author was visiting the University of Napoli ``Federico II'' and the University of Napoli Parthenope. The hospitality of these Institutions, in particular of Prof. Antonia Passarelli di Napoli and Prof. Raffaella Giova are warmly acknowledged.}

\begin{abstract}
We prove Lipschitz continuity results for solutions to a class of obstacle problems under standard growth conditions of $p-$type, $p \ge 2$. The main novelty is the use of a linearization technique going back to \cite{F85} in order to interpret our constrained minimizer as a solution to a nonlinear elliptic equation, with a bounded right hand side. This lead us to start a Moser iteration scheme which provides the $L^{\infty}$ bound for the gradient. The application of a recent higher differentiability result \cite{EPdN1} allows us to simplify the procedure of the identification of the Radon measure  in the linearization technique employed in \cite{FM00}. To our knowdledge, this is the first result for non-automonous functionals with standard growth conditions in the direction of the Lipschitz regularity.
\end{abstract}

\maketitle

\begin{center}
\fbox{\today}
\end{center}

\section{Introduction}

\vspace{3mm}
The aim of this paper is the study of the Lipschitz continuity  of the solutions $u\in W^{1,p}(\Omega)$ to variational obstacle problems of the form
 \begin{equation}
\label{obst-def0}
\min\left\{\int_\Omega f(x,Dw): w\in \mathcal{K}_{\psi, \Psi}(\Omega)\right\}.
\end{equation}
The function $\psi:\,\Omega \rightarrow [- \infty, + \infty)$, called \textit{obstacle}, belongs to the Sobolev space $W^{1,p}(\Omega)$ and the class $\mathcal{K}_{\psi,\Psi}(\Omega)$ is defined as follows, for a given $\Psi \in W^{1,p}(\overline{\Omega})$ providing the boundary value 

\vspace{-2mm}
\begin{equation}
\label{classeA}
\mathcal{K}_{\psi, \Psi}(\Omega) := \left \{w \in W^{1,p}(\Omega): w(x) \ge \psi(x) \,\, \textnormal{a.e. in $\Omega$ and } \, w|_{\partial \Omega} = \Psi|_{\partial \Omega} \right\}.    
\end{equation} 

\vspace{2mm}
In the following, we shall assume that $\mathcal{K}_{\psi, \Psi}(\Omega)$ is not empty. 

The study of the regularity theory for obstacle problems is a classical topic in Calculus of Variations and Partial Differential Equations. The first results concerning obstacle problems date back to the sixties to the pioneering work by G. Stampacchia \cite{S64} and G. Fichera \cite{F64}. The well known fact that solutions to the obstacle problem cannot be of class $\mathcal{C}^2$ independently of the regularity of the obstacle led to the origin of the concept of {\sl weak solution} and to the theory of {\it variational inequalities}, after the fundamental work of J.L. Lions and G. Stampacchia \cite{LS67}. These problems can generally be solved by means of nowadays classical methods of functional analysis; in this respect the question to establish conditions so that weak solutions can be actually classical ones is of crucial importance (see \cite{C76}). 
For more details we refer to the following monographs \cite{BC}, \cite{C84}, \cite{F82}, \cite{KS80}, \cite{R87}. 

The regularity of solutions to the obstacle problems is influenced by the one of the obstacle: in the linear case, for instance,  obstacle and solutions share the same regularity (\cite{BK74}, \cite{CK80}, \cite{KS80}), while in the nonlinear setting the situation is more involved. This generated along the years an intense research activity aimed to establish the regularity of the solution compared with the regularity assumed for the obstacle. In particular many results concern the H\"older continuity of solutions to the obstacle problem when the obstacle itself is H\"older continuous: see for instance \cite{MZ86}, \cite{C91}, \cite{E07}, \cite{EH08}, \cite{EH11}, \cite{BCO}, \cite{O16-1}; see also \cite{CL91}, \cite{L88}, \cite{FM00}, \cite{BFM01}, \cite{EH10}, \cite{EHL13}, \cite{O16-2}, \cite{B14}, \cite{O17},  \cite{BDM}, \cite{BLS15}, \cite{EPdN1}, \cite{EPdN2}.

As far as we know, such analysis has not been carried out in the direction of obtaining Lipschitz regularity results. The reason may be found in the difficulty of choosing the proper test function, in order to start the iteration process, necessary to get the desired regularity; this nowadays classical technique  (\cite{Moser}) has been adapted along the years in several situations, also in the case of non standard growth conditions, see also for instance \cite{M1}, \cite{M2}, \cite{M3}, \cite{M4}.

In \cite{FM00} Fuchs and Mingione considered constrained local minimizers of elliptic variational integrals, assuming for the Lagrangian $f$ nearly linear growth. In this respect, they were able to prove several regularity results; the one we are mainly interested in is  the Lipschitz continuity of these constrained minimizers. This step has been achived through a linearization technique, going back to \cite{F85}, later refined in \cite{D87}, see also \cite{F90}, \cite{F94}, \cite{FG}.

The main idea of the linearization approach is the following: the constrained minimizer is interpreted as a solution to an elliptic equation with a bounded right hand side, which is obtained after the identification of a suitable Radon measure. At this point the main obstruction is overcome, and higher regularity for the solution to the obstacle problem can be shown. 

We would like to remark that in \cite{FM00} a lot of effort has been employed to identify the Radon measure and the authors explicitly say that this procedure could be significantly simplified if we would have a priori proved higher differentiability for local minimizers of the obstacle problem.

However this is exactly the case for our specific situation: in a very recent paper \cite{EPdN1} Eleuteri and Passarelli di Napoli were able to establish the higher differentiability of integer and fractional order of the solutions to a class of obstacle problems (involving $p-$harmonic operators) assuming that the gradient of the obstacle possesses an extra (integer or fractional) differentiability property.
We will make use of this result to simplify the procedure outlined in \cite{FM00} and obtain the desired result.

The paper is organized as follows: Section \ref{due} contains some notations and the statement of the main result; Section \ref{tre} contains the detailed construction of the linearization procedure while the conclusion comes with Section \ref{quattro} where the a priori estimate is presented.

We would like to mention that, in order to let the linearization technique be the main focus of our work, we assumed standard growth conditions for the lagrangians of our integral functionals and $L^{\infty}$ dependence with respect to the $x-$variable. In a similar way, we assume on the  gradient of the obstacle $W^{1, \infty}$ regularity. We expect that future improvements can be obtained both in the direction of considering non-standard growth conditions (see for instance the recent papers \cite{EMM}, \cite{CGGP1}, \cite{CGGP2}), for which a counterpart of \cite{EPdN1} would be needed in the nonstandard setting, see \cite{G19}) and in the direction of weakening the assumptions on the obstacle and on the partial map $x \mapsto D_{\xi} f(x, \xi)$. We plan to continue our research in these directions.

\section{Notations and statement of the main result}

\label{due}

For a domain $\Omega \subset \mathbb{R}^n$, $n \ge 2$, and functions $u: \Omega \rightarrow \mathbb{R},$ we consider the variational integral
\begin{equation}
\label{main-F}
\mathcal{F}(u) := \int_{\Omega} f(x, Du) \, dx,
\end{equation}
where $f: \Omega \times \mathbb{R}^{n} \rightarrow [0, + \infty)$ is a Carath\'eodory function  which is convex and of class $\mathcal{C}^2$ with respect to the second variable. Moreover, we suppose that there exist two positive constants $\lambda, \Lambda$ such that, for all $\mu, \xi \in \mathbb{R}^{n}$, $\mu = \mu_i,$ $\xi = \xi_i$, $i = 1, 2, \dots, n$, a.e. in $\Omega$ and for $p \ge 2$
\begin{eqnarray}
\lambda \, (1 + |\xi|^2)^{\frac{p-2}{2}} |\mu|^2  \le \sum_{i,j} f_{\xi_i \xi_j}(x, \xi) \mu_i \mu_j,  
\label{HpCp} \\
|f_{\xi_i \xi_j}(x, \xi)| \le \, \Lambda \, (1 +|\xi|^2)^{\frac{p-2}{2}},  
\label{HpCq}\\[2.5mm]
\,\,|f_{\xi x}(x, \xi)| \le \, \Lambda \, (1 +|\xi|^{2})^{\frac{p-1}{2}}. 
\label{gtx}
\end{eqnarray}
Throughout the paper we will denote by $B_{\rho}$ and $B_R$ balls of radii respectively $\rho$ and $R$ (with $\rho < R$) compactly contained in $\Omega$ and with the same center, let us say $x_0 \in \Omega$.

Our main result reads as follows
\begin{theorem}
\label{main-result}
Let $u \in \mathcal{K}_{\psi, \Psi}(\Omega)$ be a solution to the obstacle problem \eqref{obst-def0}, under the assumptions \eqref{HpCp}-\eqref{gtx}. 
If $\psi \in W^{2, \infty}_{\rm loc}(\Omega)$, then $u \in W^{1, \infty}_{\loc}(\Omega)$ and the following estimate 
\begin{equation}
\label{a-priori}
\|(1+|Du|^2)^{\frac{1}{2}}\|_{L^{\infty}(B_{\rho})} \le \, \frac{C}{(R-\rho)^{\beta}} \, \left ( \int_{B_{R}} (1 + |Du|^2)^{\frac{p}{2}} \, dx \right )^{\frac{1}{p}} 
\end{equation}
holds for every $0 < \rho < R$ and with positive constants $C$ and $\beta$ depending on $n, p, \lambda, \Lambda, R, \rho$ and the local bounds for $\|D \psi\|_{W^{1, \infty}}$. 
\end{theorem}

\begin{remark}
Due to the local nature of our main result, we do not assume further regularity on the boundary datum $\Psi$.
\end{remark}

\section{Proof of Theorem \ref{main-result}:  the linearization procedure}

\label{tre}

Consider a smooth function $h_{\varepsilon} : (0, \infty) \to [0, 1]$ such that $h'(s) \le 0 $ for all $s\in (0, \infty)$ and 
\[   
h(s) = 
     \begin{cases}
       1 &\quad\text{for} \quad s \le \varepsilon ,\\
       0 &\quad\text{for} \quad s \ge 2\varepsilon .\\
     \end{cases}
\]

We remark that $u\in W^{1,p}(\Omega)$ is a solution to the obstacle problem in $\mathcal{K}_{\psi, \Psi}(\Omega)$ if and only if $u \in \mathcal{K}_{\psi, \Psi} $ solves the following variational inequality
\begin{equation}
\label{obst-def}
\int_{\Omega}  D_{\xi} f(x, Du) \cdot D(\varphi - u) \, dx \ge 0,
\end{equation}
for all $\varphi \in \mathcal{K}_{\psi,\Psi}$. In particular the function 
\begin{equation*}
    \varphi = u+t\cdot \eta \cdot  h_{\varepsilon}(u-\psi)  
\end{equation*}
with $\eta \in C_{0}^1 (\Omega) $, $\eta \ge 0$ and $0 < t <<1$, is such that $\varphi \in \mathcal{K}_{\psi,\Psi}$. Using $\varphi$ as test function, in the variational inequality \eqref{obst-def} gives
\begin{equation*}
     \int_{\Omega} D_{\xi} f(x, Du) \cdot D(\eta h_{\varepsilon}(u-\psi)) \, dx \ge 0  \qquad \forall  \, \eta \in C_{0}^1(\Omega) \, .
\end{equation*}
Since 
\[
\eta \mapsto L(\eta)=\int_{\Omega} D_{\xi} f(x, Du) \cdot D(\eta h_{\varepsilon}(u-\psi)) \, dx
\] 
is a bounded positive linear functional, by the Riesz representation theorem there exists a nonnegative measure $\lambda_{\varepsilon}$ such that 
\begin{equation*}
     \int_{\Omega} D_{\xi} f(x, Du) \cdot D(\eta h_{\varepsilon}(u-\psi)) \, dx = \int_{\Omega} \eta d\lambda_{\varepsilon}  \qquad \forall  \, \eta \in C_{0}^1(\Omega) \, .
\end{equation*}
We show that the measure $\lambda_{\varepsilon}$ is independent from $\varepsilon$: indeed let $\varepsilon < \varepsilon'$. The variation
\begin{equation*}
    \varphi_{+}=u+t\cdot \eta \cdot (h_{\varepsilon}-h_{\varepsilon'})(u-\psi) 
\end{equation*}
is admissibile for all $\eta \in C_{0}^1(\Omega)$ and $0<t<<1$. This gives
\begin{equation*}
    \int_{\Omega} D_{\xi} f(x, Du) \cdot D(\eta (h_{\varepsilon}-h_{\varepsilon'})(u-\psi)) \, dx \ge 0  \qquad \forall  \, \eta \in C_{0}^1(\Omega) \, .
\end{equation*}
It's easy to see that $(h_{\varepsilon}-h_{\varepsilon'})(s)$ is supported on $[\varepsilon, 2\varepsilon']$, which implies
\begin{equation*}
    \int_{\{\varepsilon < u-\psi < 2\varepsilon' \}} D_{\xi} f(x, Du) \cdot D(\eta (h_{\varepsilon}-h_{\varepsilon'})(u-\psi)) \, dx \ge 0  \qquad \forall  \, \eta \in C_{0}^1(\Omega) \, .
\end{equation*}
Then, for $t<\frac{\varepsilon}{2},$ also the variation 
\begin{equation*}
    \varphi_{-}=u-t\cdot \eta \cdot (h_{\varepsilon}-h_{\varepsilon'})(u-\psi) 
\end{equation*}
is an admissibile test function; by the same argument we get
\begin{equation*}
    \int_{\{\varepsilon < u-\psi < 2\varepsilon' \}} D_{\xi} f(x, Du) \cdot D(\eta (h_{\varepsilon}-h_{\varepsilon'})(u-\psi)) \, dx \le 0  \qquad \forall  \, \eta \in C_{0}^1(\Omega) \, ;
\end{equation*}
 together with the previous result, this implies 
\begin{equation*}
    \int_{\Omega} D_{\xi} f(x, Du) \cdot D(\eta (h_{\varepsilon}-h_{\varepsilon'})(u-\psi)) \, dx = 0  \qquad \forall  \, \eta \in C_{0}^1(\Omega) \, ,
\end{equation*}
that is equivalent to 
\begin{equation}
\label{uguale}
    \int_{\Omega} \eta d\lambda_{\varepsilon}=\int_{\Omega} \eta d\lambda_{\varepsilon'} \qquad \forall  \, \eta \in C_{0}^1(\Omega) \, .
\end{equation}
Now we can obtain equality between measures by a standard approximation argument, let $\phi \in C_{0}^{\infty}(\Omega)$ be a smooth kernel with $\phi \ge 0 $ and $\int_{\Omega} \phi \, dx = 1$, consider the corresponding family of mollifiers $(\phi_s)_{s>0}$ and any compact set $K$ such that $K \subset \subset \Omega$. For every $s>0$ such that $s<\text{dist}(K, \Omega)$ we take in equation \eqref{uguale} 
\begin{equation*}
    \eta = \chi_{K} *\phi_s  
\end{equation*}
letting $ s \searrow 0 $ gives
\begin{equation*}
    \int_{K} d\lambda_{\varepsilon} = \int_{K} d\lambda_{\varepsilon'} \, , 
\end{equation*}
in conclusion $\lambda_{\varepsilon} = \lambda_{\varepsilon'}$. \\
At this point we can rewrite our representation equation without the $\varepsilon$ dependence on the measure
\begin{equation*}
    \int_{\Omega} D_{\xi} f(x, Du) \cdot D(\eta h_{\varepsilon}(u-\psi)) \, dx = \int_{\Omega} \eta \, d\lambda  \qquad \forall  \, \eta \in C_{0}^1(\Omega) \, .
\end{equation*}
By Theorem 2.1 in \cite{EPdN1}, that holds in particular under the assumptions \eqref{HpCp}-\eqref{gtx}, we have that
\[
V_p(Du) := (1 + |Du|^2)^{\frac{p-2}{2}} Du \in W^{1,2}_{\rm loc}(\Omega),
\]
therefore we can integrate by parts and get 
\begin{equation*}
    \int_{\Omega} -\text{div}(D_{\xi} f(x, Du) ) \eta h_{\varepsilon}(u-\psi) \, dx = \int_{\Omega} \eta \, d\lambda \qquad \forall  \, \eta \in C_{0}^1(\Omega) \, .
\end{equation*}
Now, in order to identify the measure $\lambda$, we may pass to the limit  as $\varepsilon \searrow 0 $
\begin{equation}
\label{eq-g}
    \int_{\Omega} -\text{div}(D_{\xi} f(x, Du) ) \chi_{ \left[u = \psi \right] } \eta \, dx = \int_{\Omega} \eta \, d\lambda \qquad \forall  \, \eta \in C_{0}^1(\Omega) \, .
\end{equation} 
 Let us introduce 
 \begin{equation}
 \label{def-g}
 g:=\text{div}(D_{\xi} f(x, Du) ) \chi_{ \left[u = \psi \right] };
 \end{equation}
 then the equation above takes the form:
\begin{equation*}
    -\int_{\Omega} g \eta \, dx = \int_{\Omega} \eta \, d\lambda \qquad \forall  \, \eta \in C_{0}^1(\Omega) \, .
\end{equation*}
Now observing that 
\begin{equation*}
    \int_{\Omega} D_{\xi} f(x, Du) \cdot D(\eta (1-h_{\varepsilon})(u-\psi)) \, dx = 0  \qquad \forall  \, \eta \in C_{0}^1(\Omega) \, 
\end{equation*}
since $(1-h_{\varepsilon})(s)$ has support $[\varepsilon, +\infty)$, combining our results we get 
\begin{equation}
\label{eq:g}
    \int_{\Omega} D_{\xi} f(x, Du) \cdot D\eta \, dx = -\int_{\Omega} g \eta \, dx \qquad \forall  \, \eta \in C_{0}^1(\Omega) \, .
\end{equation}

We are left to obtain an $L^{\infty}$ estimate for $g$: since $Du=D\psi$ a.e. on the contact set, by (2.3) and (2.4) and the assumption $D\psi \in W^{1,\infty}_{\rm loc}(\Omega; \mathbb{R}^n)$,  we have 
\begin{equation*}
\begin{split}
    |g| &= \left|\text{div}(D_{\xi} f(x, Du) ) \chi_{ \left[u = \psi \right] } \right| =   \left|\text{div}(D_{\xi} f(x, D\psi) ) \right| \\ 
    &\le \sum_{k=1}^n | f_{\xi_k x_k}(x, D\psi) | + \sum_{k,i=1}^n | f_{\xi_k \xi_i}(x, D\psi) \psi_{x_k x_i} | \\ 
    &\le \Lambda(1+|D\psi|^2)^{\frac{p-1}{2}}+ \Lambda(1+|D\psi|^2)^{\frac{p-2}{2}}|D^2 \psi| \, 
\end{split}
\end{equation*}
that is $g \in L^{\infty}_{\rm loc}(\Omega)$.

\section{Proof of Theorem \ref{main-result}: a priori estimate}

\label{quattro}

Our starting point is now \eqref{eq:g}. We make use of the supplementary assumption $u \in W^{1, \infty}_{\rm loc}(\Omega)$, which is needed in order to let \eqref{eq:g} to be satisfied; this assumption will be removed by means of the approximation procedure. By the standard techniques of the difference quotients (see for instance \cite{giaquinta}, \cite{giusti}) we get 
\begin{equation}
\label{regE}
u \in  W^{2,2}_{\rm loc}(\Omega) \cap (1 + |Du|^2)^{\frac{p-2}{2}} |D^2 u|^2 \in L^1_{\rm loc}(\Omega)
\end{equation}
so that the ``second variation'' system holds
\begin{equation}
\label{second_variation}
\int_{\Omega} \left (\sum_{i,j=1}^n  f_{\xi_i \xi_j}(x, Du)u_{x_j x_s}D_{x_i} \varphi  \, dx + \sum_{i=1}^n f_{\xi_i x_s}(x, Du) D_{x_i} \varphi \right ) = \int_{\Omega} g \, D_{x_s} \varphi \, dx,
\end{equation}
for all $s = 1, \dots, n$ and for all $\varphi \in W^{1,p}_0(\Omega)$; here $g$ is the function which has been introduced in \eqref{def-g}. We fix $0 < \rho < R$ with $B_R$ compactly contained in $\Omega$ and we choose $\eta \in \mathcal{C}^1_0(\Omega)$ such that $0 \le \eta \le 1$,  $\eta \equiv 1$ on $B_{\rho}$, $\eta \equiv 0$ outside $B_R$ and $|D \eta| \le \frac{C}{(R - \rho)}.$ We test \eqref{second_variation} with $\varphi=\eta^2 (1+|Du|^2)^{\gamma} u_{x_s}$, for some $\gamma \ge 0$ (which is possible due to our a priori assumption \eqref{regE}) so that 
$$ D_{x_i} \varphi =2 \eta \eta_{x_i} (1 +|Du|^2)^{\gamma} u_{x_s} + 2 \eta^2 \gamma (1 +|Du|^2)^{\gamma-1} |Du| D_{x_i}(|Du|) u_{x_s} + \eta^2 (1 +|Du|^2)^{\gamma} u_{x_s x_i}$$
Inserting in \eqref{second_variation} we get: 
\begin{eqnarray*} 
0 &=& \int_{\Omega} \sum_{i,j=1}^n f_{\xi_i \xi_j}(x, Du) u_{x_j x_s} 2\eta \eta_{x_i} (1 +|Du|^2)^{\gamma} u_{x_s} \, dx\\ 
&& +  \int_{\Omega} \sum_{i,j=1}^n f_{\xi_i \xi_j}(x, Du) u_{x_j x_s} \eta^2 (1 +|Du|^2)^{\gamma} u_{x_s x_i} \, dx\\
&& +  \int_{\Omega} \sum_{i,j=1}^n f_{\xi_i \xi_j}(x, Du) u_{x_j x_s} 2\eta^2 \gamma  (1 +|Du|^2)^{\gamma-1} |Du| D_{x_i}(|Du|) u_{x_s} \, dx\\
&& + \int_{\Omega} \sum_{i=1}^n f_{\xi_i x_s}(x, Du) 2\eta \eta_{x_i}  (1 +|Du|^2)^{\gamma} u_{x_s} \, dx\\
&& +  \int_{\Omega} \sum_{i=1}^n f_{\xi_i x_s }(x, Du) \eta^2 (1 + |Du|^2)^{\gamma}  u_{x_s x_i} \, dx\\
&& +  \int_{\Omega} \sum_{i=1}^n f_{\xi_i x_s}(x, Du) 2\eta^2 \gamma  (1 + |Du|^2)^{\gamma-1}  |Du| D_{x_i}(|Du|) u_{x_s} \, dx \\
&& - \int_{\Omega} g  2 \eta  \eta_{x_s} \,  (1 + |Du|^2)^{\gamma}  \, u_{x_s} \, dx\\
&& - \int_{\Omega} g  2 \eta^2 \gamma \,  (1 + |Du|^2)^{\gamma - 1}  |Du| \, D_{x_s}(|Du|) \, u_{x_s} \, dx \\
&& - \int_{\Omega} g \, \eta^2  (1 + |Du|^2)^{\gamma} u_{x_s x_s} \, dx\\
&& =: I_{1,s} + I_{2,s} + I_{3,s} + I_{4,s} + I_{5,s} + I_{6,s} + I_{7,s} + I_{8,s} + I_{9,s}.
\end{eqnarray*}
We sum in the previous  equation  all terms with respect to $s$ from 1 to $n$, and we denote by $I_1-I_9$ the corresponding integrals. In the sequel constants will be denoted by $C$, regardless their actual value. Only the relevant dependencies will be highlighted. 
\\
By the Cauchy-Schwarz inequality, the Young inequality and \eqref{HpCq}, we have
\begin{eqnarray*}
|{I}_1| &=& \left |\int_{\Omega} 2 \eta  (1 + |Du|^2)^{\gamma}   \sum_{i,j,s =1}^n f_{\xi_i \xi_j}(x, Du) u_{x_j x_s} \eta_{x_i} u_{x_s} \, dx\right |\\
& = &
\Bigg| \int_{\Omega} 2 \eta  (1 + |Du|^2)^{\gamma} \\
&& \times \sum_{s=1}^n \left( \sum_{i,j=1}^n f_{\xi_i \xi_j}(x, Du) \eta_{x_i} \eta_{x_j} u_{x_s}^2 \right)^{1/2} \left(\sum_{i,j=1}^n 
f_{\xi_i \xi_j}(x, Du) u_{x_s x_i} \, u_{x_s x_j} \right)^{1/2} \, dx \Bigg| \\ 
&\le& \int_{\Omega} 2 \eta  (1 + |Du|^2)^{\gamma} \\
&& \times \left \{ \sum_{i,j,s =1}^n f_{\xi_i \xi_j}(x, Du) \eta_{x_i} \eta_{x_j} u_{x_s}^2\right \}^{1/2} \, \left \{ \sum_{i,j,s =1}^n f_{\xi_i \xi_j}(x, Du) u_{x_s x_i} \, u_{x_s x_j}\right \}^{1/2} \, dx \\
&\le &  C   \int_{\Omega}  (1 + |Du|^2)^{\gamma}  \sum_{i,j,s =1}^n f_{\xi_i \xi_j}(x, Du) \eta_{x_i} \eta_{x_j} u_{x_s}^2 \, dx\\
&&  + \frac{1}{4} \int_{\Omega} \eta^2  (1 + |Du|^2)^{\gamma}  \sum_{i,j,s =1}^n f_{\xi_i \xi_j}(x, Du) u_{x_s x_i} \, u_{x_s x_j}  \, dx \\
&\le &  C (\Lambda)  \int_{\Omega} \Lambda  (1 + |Du|^2)^{\frac{p-2}{2}+\gamma}  \sum_{i,j,s =1}^n \eta_{x_i} \eta_{x_j} u_{x_s}^2 \, dx \\
&& + \frac{1}{4} \int_{\Omega} \eta^2  (1 + |Du|^2)^{\gamma}  \sum_{i,j,s =1}^n f_{\xi_i \xi_j}(x, Du) u_{x_s x_i} \, u_{x_s x_j} \, dx \\
&\le &  C(n, \Lambda)  \int_{\Omega} |D\eta|^2  (1 + |Du|^2)^{\frac{p-2}{2}+\gamma}  \sum_{s =1}^n u_{x_s}^2 \, dx \\
&& + \frac{1}{4} \int_{\Omega} \eta^2  (1 + |Du|^2)^{\gamma}  \sum_{i,j,s =1}^n f_{\xi_i \xi_j}(x, Du) u_{x_s x_i} \, u_{x_s x_j} \, dx \\
&\le& C \int_{\Omega} |D \eta|^2 \,  (1 + |Du|^2)^{p/2+\gamma}  \, dx \\
&& + \frac{1}{4} \int_{\Omega} \eta^2  (1 + |Du|^2)^{\gamma}  \, \sum_{i,j,s =1}^n f_{\xi_i \xi_j}(x, Du) u_{x_j x_s} u_{x_i x_s} \, dx.
\end{eqnarray*}
On the other hand, using \eqref{HpCp} and the fact that $D_{x_j}(|Du|)|Du|=\sum_{k=1}^n u_{x_j x_k} u_{x_k}$, we can estimate the term $I_3$ as follows: 
\begin{eqnarray*}
|{I}_3| &=& \int_{\Omega}  \sum_{i,j,s =1}^n f_{\xi_i \xi_j}(x, Du) u_{x_j x_s} \left [2\eta^2 \gamma  (1 + |Du|^2)^{\gamma-1}  D_{x_i}(|Du|)|Du| \right ] u_{x_s} \, dx\\ 
&\ge& 2 \gamma \, \int_{\Omega} \eta^2 |Du|^{2\gamma-1} \sum_{i,j,s =1}^n f_{\xi_i \xi_j}(x, Du) D_{x_i}(|Du|) u_{x_j x_s} u_{x_s} \, dx\\
&=& 2 \gamma \, \int_{\Omega} \eta^2 |Du|^{2\gamma-1} \sum_{i,j =1}^n f_{\xi_i \xi_j}(x, Du) D_{x_i}(|Du|) \left(\sum_{s=1}^n u_{x_j x_s} u_{x_s} \right)\, dx\\
&=& 2 \gamma \, \int_{\Omega} \eta^2  |Du|^{2\gamma} \sum_{i,j =1}^n f_{\xi_i \xi_j}(x, Du) D_{x_i}(|Du|) D_{x_j}(|Du|) \, dx\\
&\ge& 2 \gamma \, \int_{\Omega} \eta^2  |Du|^{2\gamma} |D(|Du|)|^2  (1 + |Du|^2)^{\frac{p-2}{2}}  \, dx \ge 0. \\
\end{eqnarray*}
We can estimate the fourth and the fifth term by the Cauchy-Schwarz and the Young inequalities, together with \eqref{gtx}, as follows
\begin{eqnarray*}
|{I}_4| &=& 2 \int_{\Omega} \eta  (1 + |Du|^2)^{\gamma}  \sum_{i,s =1}^n f_{\xi_i x_s}(x, Du) \eta_{x_i} u_{x_s} \, dx\\ 
&\overset{\eqref{gtx}}{\le}&  2 \Lambda \, \int_{\Omega} \eta  \,  (1 + |Du|^2)^{\gamma+\frac{p-1}{2}}  \sum_{i,s=1}^n |\eta_{x_i} u_{x_s} | \, dx \\
&\le&   C (\Lambda)  \int_{\Omega} \eta |D\eta| |Du| \, (1 + |Du|^2)^{\gamma+\frac{p-1}{2}}   \, dx \\
&\le& C \int_{\Omega} (\eta^2+|D\eta|^2) (1 + |Du|^2)^{\gamma+\frac{p}{2}},
\end{eqnarray*} 
and also:
\begin{eqnarray*}
|{I}_5| &=& \left| \int_{\Omega} \eta^2 (1 + |Du|^2)^{\gamma} \sum_{i,s=1}^n f_{\xi_i x_s }(x, Du)  u_{x_s x_i} \, dx \right| \\  
&\overset{\eqref{gtx}}{\le}& \Lambda \, \int_{\Omega} \eta^2 (1 + |Du|^2)^{\gamma+\frac{p-1}{2}}  \, \sum_{i,s=1}^n  u_{x_s x_i} \, dx \\ 
&=& \Lambda \, \int_{\Omega} \eta^2 (1 + |Du|^2)^{\gamma+\frac{p-1}{2}}  |D^2 u| \, dx \\ 
&\le& \frac{1}{4} \int_{\Omega} \eta^2 (1 + |Du|^2)^{\frac{p-2}{2}+\gamma}  |D^2 u|^2 \, dx + C \, \int_{\Omega} \eta^2 (1 + |Du|^2)^{\frac{p}{2}+\gamma}  \, dx.
\end{eqnarray*}
Finally, by similar arguments, the Cauchy-Schwarz inequality, \eqref{HpCq} 
and $|D(|Du|)| \le |D^2 u|$ give
\begin{eqnarray*}
|{I}_6| &=&  2 \gamma \, \int_{\Omega} \sum_{i,s=1}^n f_{\xi_i x_s}(x, Du) \eta^2  (1 + |Du|^2)^{\gamma-1}  |Du| D_{x_i}(|Du|) u_{x_s} \, dx\\ 
&=&2 \gamma \, \int_{\Omega} \eta^2 (1 + |Du|^2)^{\gamma-1}  |Du| \sum_{i,s=1}^n f_{\xi_i x_s}(x, Du)  D_{x_i}(|Du|) u_{x_s} \, dx\\
&\le& 2 \, \gamma \, \int_{\Omega} \eta^2 (1+|Du|^2)^{\gamma-\frac{1}{2}}  \sum_{i,s=1}^n f_{\xi_i x_s}(x, Du)  D_{x_i}(|Du|) u_{x_s} \, dx\\
&\le& 2 \, \gamma \, \Lambda \, \int_{\Omega} \eta^2  (1 + |Du|^2)^{\gamma-\frac{1}{2}+\frac{p-1}{2}}  \sum_{i,s=1}^n D_{x_i}(|Du|) u_{x_s} \, dx\\
&\le& 2 \, n \, \gamma \, \Lambda \, \int_{\Omega} \eta^2 (1 + |Du|^2)^{\gamma-\frac{1}{2}+\frac{p-1}{2}}  |D(|Du|)| |Du| \, dx\\
&\le& 2 \, n \, \gamma \, \Lambda \, \int_{\Omega} \, \eta^2 (1 + |Du|^2)^{\gamma+\frac{p-1}{2}}  |D^2 u| \, dx\\
&\le& \frac{1}{4} \int_{\Omega} \eta^2 |D^2u|^2 (1 + |Du|^2)^{\frac{p-2}{2}+\gamma}   \, dx + C \gamma^2 \, \int_{\Omega} \eta^2  (1 + |Du|^2)^{\frac{p}{2}+\gamma }  \, dx,
\end{eqnarray*}
where the constant $C$ depends only on $n, p, \Lambda$ but it is independent of $\gamma$.
\\
Let us now deal with the terms containing the function $g$. We use the bound 
\[
\|g\|_{L^{\infty}_{\rm loc}(\Omega)} \le \, C,
\]
established in Section \ref{tre}; here the constant $C$ is independent of $n$ and $\gamma$. 
\\
We first have
\begin{eqnarray*}
|I_7| &\le& \|g\|_{L^{\infty}(B_R)} \int_{\Omega} 2 |D \eta| (1 + |Du|^2)^{\gamma}  \, |Du| \, dx \\
&\le& \, C \, \int_{\Omega} |D \eta|^2 (1 + |Du|^2)^{\gamma + p/2} \, dx.
\end{eqnarray*}
\\
Arguing similarly as in the estimate of $I_6$ we deduce
\begin{eqnarray*}
|I_8| &\le&  \, 2 \gamma \int_{\Omega} |g| \, \eta^2 (1 + |Du|^2)^{\gamma - 1}  |D(|Du|)| \, |Du|^2 \, dx \\
&\le&  \frac{1}{4} \int_{\Omega} \eta^2 |D^2u|^2 (1 + |Du|^2)^{\frac{p-2}{2}+\gamma}  \, dx + C \gamma^2 \, \|g\|^2_{L^{\infty}(B_R)} \,  \int_{\Omega} \eta^2 (1 + |Du|^2)^{\frac{p}{2}+\gamma }  \, dx.
\end{eqnarray*}
Finally
\begin{eqnarray*}
|I_9| &\le& \frac{1}{4} \int_{\Omega} \eta^2 |D^2u|^2 (1 + |Du|^2)^{\frac{p-2}{2}+\gamma}  \, dx + C  \|g\|^2_{L^{\infty}(B_R)}  \int_{\Omega} \eta^2  (1 + |Du|^2)^{\frac{p}{2}+\gamma }  \, dx.
\end{eqnarray*}
We remark that in the estimates of the terms $I_7-I_9$ we used in an essential way the fact $(1 + |Du|^2) \ge 1$. It is possible to overcome this limtation and dealing also with the degenerate case by proceeding as in \cite{EMM}. 
\\
Summing up and using \eqref{HpCp} we obtain
\[
\int_{\Omega} \, \eta^2 \, (1 + |Du|^2)^{\frac{p-2}{2} + \gamma}  |D^2 u|^2 \, dx \le \, C \, (1 + \gamma^2)  \int_{\Omega} (\eta^2 + |D \eta|^2) \, (1 + |Du|^2)^{\frac{p}{2} + \gamma}  \, dx,
\]
for any  $0 < \rho < R$, where the constant $C$ depends on $\lambda, \Lambda, n, p$ but is independent of $\gamma$.

On the other hand
\begin{eqnarray*}
&& \int_{\Omega} \left| D \left [\eta (1 + |Du|^2)^{\frac{p}{4} + \frac{\gamma}{2}}  \right ]\right|^2 \, dx \\
&\le& \int_{\Omega} |D \eta|^2 (1 + |Du|^2)^{\frac{p}{2} + \gamma} \, dx + \int_{\Omega} \eta^2 (1 + |Du|^2)^{\frac{p-2}{2} + \gamma}  |D^2 u| \, dx\\
&\le& \, C \,   (1 + \gamma^2)  \int_{\Omega} (\eta^2 + |D \eta|^2) \, (1 + |Du|^2)^{\frac{p}{2} + \gamma}  \, dx.
\end{eqnarray*}
Now, let $2^* = \frac{2n}{n - 2}$ for $n > 2$ while $2^*$ equal to any fixed real number greater than 2, if $n = 2$. By the Sobolev's inequality there exists a constant $C$ such that
\begin{eqnarray*}
\left ( \int_{\Omega} \eta^{2^{*}}  \left (1 + |Du|^2 \right )^{\left ( \frac{p}{2} + \gamma \right ) \frac{2^*}{2}} \, dx \right)^{\frac{2}{2^*}} &\le& \, C \, \int_{\Omega} \left| D \left [\eta (1 + |Du|^2)^{\frac{p}{4} + \frac{\gamma}{2}}  \right ]\right|^2 \, dx \\
&\le& \,  C \, (1 + \gamma^2)  \int_{\Omega} (\eta^2 + |D \eta|^2) \, (1 + |Du|^2)^{\frac{p}{2} + \gamma} \, dx.
\end{eqnarray*}

\section{Proof of Theorem \ref{main-result}: an iteration scheme}
By the definition of $\eta$ we obtain
\begin{equation}\label{stima main}
\left ( \int_{B_\rho} \left (1 + |Du|^2 \right )^{\left ( \frac{p}{2} + \gamma \right ) \frac{2^*}{2}} \, dx \right)^{\frac{2}{2^*}} \le  \frac{C (1 + \gamma^2)}{(R-\rho)^2}  \int_{B_R} (1 + |Du|^2)^{\frac{p}{2} + \gamma} \, dx, 
\end{equation}
for any $B_\rho$ and $B_R$ compactly contained in $\Omega$ with $\rho<R$ and the same center. \\ \noindent Let us define the sequence $(\gamma_k)_k$ by induction as follows
\begin{equation}\label{gamma def}
    \gamma_1=0 \, ; \qquad \gamma_{k+1}=\left(\gamma_k+\frac{p}{2}\right)\frac{2^{*}}{2}-\frac{p}{2}, \quad \forall \, k \ge 1 .
\end{equation}

If $(\gamma_k)_k$ is the sequence defined in \eqref{gamma def}, then the following representation formula holds
\begin{equation}\label{lemma gammak}
    \gamma_{k}=\frac{p}{2}\left(\left(\frac{2^{*}}{2}\right)^{k-1}-1 \right),
\end{equation}
the proof easily follows by induction.

Fix two rays $R_0$ and $\rho_0$ such that $0<\rho_0<R_0<\rho_0+1$, let us define $R_k=\rho_0+(R_0-\rho_0)2^{-k}$ for every $k\ge 1$. Since $R_{k+1}<R_k \,\, \forall \, k\ge 1$, we are allowed to insert $\rho=R_{k+1}$ and $R=R_k$ in estimate \eqref{stima main}, we also take $\gamma=\gamma_k$, where $(\gamma_k)_k$ has been defined in \eqref{gamma def}. With these choices, by \eqref{stima main} we obtain 
\setlength{\jot}{9pt}
\begin{align}
\left ( \int_{B_{R_{k+1}}} (1 + |Du|^2   )^{\left ( \frac{p}{2} + \gamma_k \right ) \frac{2^*}{2}} \, dx \right)^{\frac{2}{2^*}}  &\le  \frac{C (1 + \gamma_k^2)}{(R_k-R_{k+1})^2}  \int_{B_{R_k}} (1 + |Du|^2)^{\frac{p}{2} + \gamma_k} \, dx \nonumber \\ &= \frac{C(1 + \gamma_k^2)4^{k+1}}{(R_0-\rho_0)^2} \int_{B_{R_k}} (1 + |Du|^2)^{\frac{p}{2} + \gamma_k} \, dx , \label{stima main con rk}
\end{align}
for every $k\ge 1$. Since the constant $C$ was independent of $\gamma$, it is clear that now $C$ is independent on $k$. 
Let us define 
\begin{equation}
    A_k= \left(\int_{B_{R_k}} (1 + |Du|^2)^{\frac{p}{2} + \gamma_k} \, dx \right)^{\frac{1}{p+2\gamma_k}}, \quad \forall \, k\ge 1.
\end{equation}
The estimate \eqref{stima main con rk} directly gives us an inductive estimate on $(A_k)_k$. By definition \eqref{gamma def} of $\gamma_k$, and by \eqref{stima main con rk} we have 
\setlength{\jot}{8pt}
\begin{align}
    A_{k+1} &= \left(\int_{B_{R_{k+1}}} (1 + |Du|^2)^{\frac{p}{2} + \gamma_{k+1}} \, dx \right)^{\frac{1}{p+2\gamma_{k+1}}} \nonumber\\&= \left(\int_{B_{R_{k+1}}} (1 + |Du|^2)^{\left(\frac{p}{2} + \gamma_{k}\right)\frac{2^{*}}{2}} \, dx \right)^{\frac{1}{\left(\frac{p}{2} + \gamma_{k}\right)2^*}} \nonumber\\ &\le \left[\frac{C(1 + \gamma_k^2)4^{k+1}}{(R_0-\rho_0)^2} \right]^{\frac{1}{p+2\gamma_k}} \cdot \left(\int_{B_{R_k}} (1 + |Du|^2)^{\frac{p}{2} + \gamma_k} \, dx \right)^{\frac{1}{p+2\gamma_k}} \nonumber\\ &= \left[\frac{C(1 + \gamma_k^2)4^{k+1}}{(R_0-\rho_0)^2} \right]^{\frac{1}{p+2\gamma_k}} \cdot A_k \, , \label{stima su ak}
\end{align}
for every $k\ge 1$. By iterating \eqref{stima su ak} we obtain
\begin{equation}\label{produttoria}
    A_{k+1} \le A_1 \cdot \prod_{i=1}^k \left[\frac{C(1 + \gamma_i^2)4^{i+1}}{(R_0-\rho_0)^2} \right]^{\frac{1}{p+2\gamma_i}} \le C \cdot A_1 \cdot (R_0-\rho_0)^{-\sum_{i=1}^{\infty} 1/(p/2+\gamma_i)} , 
\end{equation}
for every $k\ge 1$. Here we used the fact that $\prod_{i=1}^k \big[C(1+\gamma_i^2)4^{i+1} \big]^{\frac{1}{p+2\gamma_i}}$ is bounded, indeed
\begin{equation*}
    \prod_{i=1}^{\infty} \big[C(1+\gamma_i^2)4^{i+1} \big]^{\frac{1}{p+2\gamma_i}} = \textnormal{exp}\left(\sum_{i=1}^{\infty} \frac{\ln(C(1+\gamma_i^2)4^{i+1})}{p+2\gamma_i} \right) < + \infty \,, 
\end{equation*}
since the series is convergent by formula \eqref{lemma gammak}.
Since $\rho_0 < R_k < R_0 $, by exploting $A_1$ in \eqref{produttoria}, for every $k\ge 2$ we deduce that
\begin{equation}\label{fine iterazione beta}
    \left(\int_{B_{\rho_0}} (1 + |Du|^2)^{\frac{1}{2}\left(p + 2\gamma_k\right)} \, dx \right)^{\frac{1}{p+2\gamma_k}} \le \frac{C}{(R_0-\rho_0)^\beta} \left(\int_{B_{R_0}} (1 + |Du|^2)^{\frac{p}{2}} \, dx \right)^{\frac{1}{p}} , 
\end{equation}
where we defined
\begin{equation*}
    \beta := \sum_{i=1}^{\infty} \frac{1}{\frac{p}{2}+\gamma_i} = \frac{2^{*}}{p(2^{*}/2)-p} .
\end{equation*}
Since $p+2\gamma_k \to +\infty $, as $k \to +\infty $ it is clear that the left hand side of \eqref{fine iterazione beta} converges to the essential supremum of $(1+|Du|^2)^{\frac{1}{2}}$ in $B_{\rho_0}$. Thus, provided the further regulatity that we will remove in the next section with an approximation procedure, we conclude the proof of Theorem \ref{main-result}. 

\section{Proof of Theorem \ref{main-result}: approximation}

First of all we state an approximation theorem for $f$ through a suitable sequence of regular functions. 
Let $B$ be the unit ball of $\mathbb{R}^n$ centered in the origin and consider a positive decreasing sequence $\varepsilon_{\ell} \rightarrow 0$. Let $\phi \in C_{0}^{\infty}(B)$ be a smooth symmetric kernel with $\phi \ge 0$ and $\int_{B} \phi \, dx = 1$. \\ \noindent Let us define 
\begin{equation}
    f^{\ell}(x, \xi) = \int_{B \times B} \phi(y) \phi(\eta) f(x + \varepsilon_{\ell} y, \xi + \varepsilon_{\ell} \eta) \, d \eta \, d y,
\end{equation}
and consequently 
\begin{equation}\label{defdifellk}
f^{\ell k}(x, \xi) = f^{\ell}(x, \xi) + \frac{1}{k} (1 + |\xi|^2)^{\frac{p}{2}}.
\end{equation}
The proof of the following theorem follows from similar arguments as in \cite{EMM}.
\begin{theorem}\label{app teorema}
Let $f$ satisfy the growth conditions \eqref{HpCp}, \eqref{HpCq}, \eqref{gtx}, and let $f(x, \xi)$ be a convex $C^2$ function in its second argument. Then the sequence $f^{\ell k}: \Omega \times \mathbb{R}^n \rightarrow [0, + \infty)$ defined in \eqref{defdifellk}, convex in the last variable is such that $ f^{\ell k} \to f $ pointwise as $ (\ell, k) \to + \infty $ for a.e. $x \in \Omega$ and for all $\xi \in \mathbb{R}^n$. Moreover, $f^{\ell k} \to f $ uniformly in $ \Omega_0 \times K $ for every $\Omega_0 \subset \subset \Omega$ and $K$ compact set of $\mathbb{R}^n$. Moreover:
\begin{itemize}
  \item There exist constants $M_1, M_2 >0 $, independent of $k,\ell$ such that 
  \begin{equation}\label{crescita p approssimate}
M_1(1+|\xi|^2)^{\frac{p}{2}} \le f^{\ell k}(x, \xi) \le M_2 (1 + |\xi|^2)^{\frac{p}{2}} \qquad \textnormal{for a.e. $x \in \Omega$, for all $\xi \in \mathbb{R}^n$},
\end{equation}

  \item there exist constants $\lambda_1, \Lambda_1>0$ such that for all $(x, \xi) \in \Omega \times \mathbb{R}^n$ and $\mu \in \mathbb{R}^n $
  \begin{equation}
  \lambda_1 \, (1 + |\xi|^2)^{\frac{p-2}{2}} |\mu|^2  \le \sum_{i,j} f^{\ell k}_{\xi_i \xi_j}(x, \xi) \mu_i \mu_j ,
  \end{equation}
  and 
  \begin{equation}
  |f^{\ell k}_{\xi_i \xi_j}(x, \xi)| \le  \Lambda_1 \, (1 +|\xi|^2)^{\frac{p-2}{2}}.
  \end{equation}
  
  \item There exists a constant $\Lambda_2>0$ such that for all $(x, \xi) \in \Omega \times \mathbb{R}^n$ 
  \begin{equation}
      |f^{\ell k}_{\xi x}(x, \xi)| \le  \Lambda_2  (1 +|\xi|^{2})^{\frac{p-1}{2}}. 
  \end{equation}
\end{itemize}
\end{theorem}
\noindent Let $u \in W^{1,p}(\Omega)$ denote a local minimizer of the problem 
\begin{equation}\label{non app problem}
\min\left\{\int_\Omega f(x,Dw): w\in \mathcal{K}_{\psi, \Psi}(\Omega)\right\}.
\end{equation}
Fix $x_0\in \Omega $ and a radius $R>0$ such that $B_{R}(x_0)\subset \subset \Omega$. Consider the following variational problem 
\begin{equation}\label{app problem}
\inf \left \{ \int_{B_R} f^{\ell k}(x, Dv) \, dx, \,\,\, v \in W^{1,p}_0(B_R) + u \, , \, v\ge \psi \,\,\, \textnormal{a.e. in $\Omega$} \, \right\},
\end{equation} 
where $f^{\ell k}$ have been defined in \eqref{defdifellk}. By classical lower semicontinuity arguments, there exist $v^{\ell k} \in W^{1,p}_0(B_R) + u $ solution to problem \eqref{app problem}. By the growth conditions \eqref{crescita p approssimate} and the minimality of $v^{\ell k}$, we get
\setlength{\jot}{10pt}
\begin{align*}
\int_{B_R} |D v^{\ell k}|^p \, dx &\le \int_{B_R} (1+|D v^{\ell k}|^2)^{\frac{p}{2}} \, dx \\ &\le C \int_{B_R} f^{\ell k}(x, Dv^{\ell k}) \, dx \le  C \int_{B_R} f^{\ell k}(x, D u) \, dx \\
&= C \int_{B_R} f^{\ell} (x, D u) \, dx + \frac{C}{k} \int_{B_R} (1 + |Du|^2)^{\frac{p}{2}} \, dx \\ &\le C \int_{B_R} (1+|Du|^2)^{\frac{p}{2}} \, dx .
\end{align*}
Note that the right hand side of the previous estimate does not depend on $\ell$, thus $\|Dv^{\ell k}\|_{L^p(B_R)}$ is bounded for every fixed $k$. Since $W^{1,p}(B_R)$ is weakly compact, up to passing to a subsequence, we can assume there exist $v^{k}\in u+W^{1,p}_0(B_R)$ such that 
\begin{equation}
    v^{\ell k} \stackrel{\ell \rightarrow \infty}{\rightarrow} v^k \,\,\, \textnormal{weakly in $W^{1,p}_0(B_R) + u$},
\end{equation}
where, for every $k\ge 1$, the limit function $v^k$ still belongs to $\mathcal{K}_{\psi, \Psi}$ since this set is weakly closed. 
Moreover, by classical convolution properties $f^{\ell}(x, Du) \to f(x, Du) $ pointwise a.e. in $B_R$ as $\ell \to +\infty $, and again by growth conditions \eqref{crescita p approssimate} 
\begin{equation}\label{domination}
    \int_{B_R} f^{\ell}(x, Du) \, dx \le  C \int_{B_R} (1 + |D u|^2)^{\frac{p}{2}} \, dx,
\end{equation}
for every $\ell \ge 1$. Thus by the Dominated Convergence Theorem we deduce 
\begin{equation}\label{domination estimate}
    \lim_{\ell \rightarrow \infty} \int_{B_R} |D v^{\ell k}|^p \, dx 
\le C \int_{B_R} f(x, D u) \, dx + \frac{C}{k} \int_{B_R} (1 + |Du|^2)^{\frac{p}{2}} \, dx.
\end{equation}
Now, in force of Theorem \ref{app teorema}, we know that $f^{\ell k}$ satisfies the growth conditions needed to apply our a priori estimate for every $\ell, k$. Thus, applying the a priori estimate to the solutions to \eqref{app problem}, using the minimality of $v^{\ell k}$
\setlength{\jot}{10pt}
\begin{align*}
    \|Dv^{\ell k}\|_{L^{\infty}(B_\rho)} & \le \widetilde{C} \left( \int_{B_R} (1+|Dv^{\ell k}|^2)^{\frac{p}{2}} \, dx \right)^{\frac{1}{p}} \\ &\le \widetilde{C} \left( \int_{B_R} f^{\ell k}(x, Dv^{\ell k}) \, dx \right)^{\frac{1}{p}}  \\ &\le \widetilde{C} \left( \int_{B_R} f^{\ell k}(x, Du) \, dx \right)^{\frac{1}{p}} \\ &\le \widetilde{C} \left( \int_{B_R} (1+|Du|^2)^{\frac{p}{2}} \, dx \right)^{\frac{1}{p}} ,
\end{align*}
where the constant $\widetilde{C}=\frac{C}{(R-\rho)^{\beta}}$ is independent of $\ell, k$. Therefore, for every fixed $k$, the sequence $\|Dv^{\ell k}\|_{L^{\infty}(B_R)}$ is bounded. Since weak limits are unique, we conclude that 
\begin{equation}
    v^{\ell k} \stackrel{\ell \rightarrow \infty}{\rightarrow} v^k \,\,\, \textnormal{weakly star in $W^{1,\infty}_{\rm loc}(B_R)$}.
\end{equation}
By estimate \eqref{domination estimate} 
\setlength{\jot}{10pt}
\begin{align}
    \|Dv^{k}\|_{L^{p}(B_R)} &\le \liminf_{\ell \to \infty} \|Dv^{\ell k}\|_{L^{p}(B_R)} \nonumber \\ &\le C \left( \int_{B_R} f(x, D u) \, dx + \frac{1}{k} \int_{B_R} (1 + |Du|^2)^{\frac{p}{2}} \, dx \right)^{\frac{1}{p}} \nonumber \\ &\le C \left( \int_{B_R} f(x, D u)+ (1 + |Du|^2)^{\frac{p}{2}} \, dx \right)^{\frac{1}{p}} , 
\end{align}
for every $k \ge 1$. Moreover, for any $\rho$ such that $0<\rho <R$, we have
\setlength{\jot}{10pt}
\begin{align}
    \|Dv^{k}\|_{L^{\infty}(B_\rho)} &\le \liminf_{\ell \to \infty} \|Dv^{\ell k}\|_{L^{\infty}(B_\rho)} \nonumber \\ &\le \widetilde{C} \left( \int_{B_R} (1+|Du|^2)^{\frac{p}{2}} \, dx \right)^{\frac{1}{p}} . \label{stima per pippo}
\end{align}
Thus we can deduce that there exists $\bar{v} \in u + W^{1,p}_0(B_R)$ such that, up to subsequences
\begin{equation*}
     v^{k} \rightarrow \bar{v} \,\,\, \textnormal{weakly in $W^{1,p}_0(B_R) + u$},\qquad \qquad v^{k} \rightarrow \bar{v} \,\,\, \textnormal{weakly star in $W^{1,\infty}_{\rm loc}(B_R)$}
\end{equation*}
and the limit function $\bar{v}$ still belongs to $\mathcal{K}_{\psi, \Psi}$ since this set is weakly closed. Since estimate \eqref{stima per pippo} holds for every $k\ge 1$, for any $0< \rho < R$ we obtain 
\begin{equation}
   \|D\bar{v}\|_{L^{\infty}(B_\rho)} \le \frac{C}{(R-\rho)^\beta} \left( \int_{B_R} (1+|Du|^2)^{\frac{p}{2}} \, dx \right)^{\frac{1}{p}} .
\end{equation}

\noindent Hence combining lower semicontinuity, the minimality of $v^{\ell k}$ for problem \eqref{app problem}, and properties of mollification we obtain
\setlength{\jot}{10pt}
\begin{align*}
\int_{B_{\rho}} f(x, D \bar{v}) \, dx &\le \liminf_{k \to \infty} \int_{B_{\rho}} f(x, D v^{k}) \, dx \\&\le \liminf_{k \to \infty} \liminf_{\ell \to \infty} \int_{B_{\rho}} f(x, D v^{\ell k}) \, dx \\ &= \liminf_{k \to \infty} \liminf_{\ell \to \infty} \int_{B_{\rho}} f^{\ell}(x, D v^{\ell k}) \, dx \\ &\le  \liminf_{k \to \infty} \liminf_{\ell \to \infty} \left(\int_{B_R} f^{\ell}(x, D v^{\ell k})+\frac{1}{k}(1+|Dv^{\ell k}|^2)^{\frac{p}{2}} \, dx \right) \\&= \liminf_{k \to \infty} \liminf_{\ell \to \infty} \int_{B_{R}} f^{\ell k}(x, D v^{\ell k}) \, dx \\ &\le \liminf_{k \to \infty} \liminf_{\ell \to \infty} \int_{B_{R}} f^{\ell k}(x, D u) \, dx \\ &= \liminf_{k \to \infty} \liminf_{\ell \to \infty} \left(\int_{B_R} f^{\ell}(x, D u) \, dx +\frac{1}{k} \int_{B_R} (1+|Du|^2)^{\frac{p}{2}} \, dx \right) \\ &= \liminf_{k \to \infty}  \left(\int_{B_R} f(x, D u) \, dx +\frac{1}{k} \int_{B_R} (1+|Du|^2)^{\frac{p}{2}} \, dx \right) \\ &= \int_{B_R} f(x, D u) \, dx , 
\end{align*}
for every $0<\rho<R$. Then, letting $\rho \rightarrow R$ in the previous inequality, we get 
\begin{equation}
    \int_{B_R} f(x, D \bar{v}) \, dx \le \int_{B_R} f(x, D u) \, dx .
\end{equation}
Therefore, $u$ and $\bar{v}$ are two solutions to problem \eqref{non app problem}. We claim that $u=\bar{v}$ on $B_R$. Indeed, if we suppose that $u \neq {v}$ choose any $\theta \in (0,1)$ and define $v_{\theta}=\theta \bar{v} + (1-\theta) u $, since $f(x, \xi)$ is strictly convex in its second argument we have 
\begin{align*}
    \int_{B_R} f(x, Dv_{\theta}) \, dx &< \theta \int_{B_R} f(x, D\bar{v}) \, dx + (1-\theta) \int_{B_R} f(x, Du) \, dx \\ &= \int_{B_R} f(x, Du) \, dx ,
\end{align*}
contradicting the minimality of $u$. Thus, regularity is preserved through the limit procedure, and this concludes the proof of Theorem \ref{main-result} because we have removed the further regularity assumptions we needed to obtain the result.


\begin{thebibliography}{99}

\bibitem{BC} \textsc{C. Baiocchi, A. Capelo:} \emph{Disequazioni variazionali e quasi variazionali. Applicazioni a problemi di frontiera libera,} Quaderni U.M.I. Pitagora, Bologna 1978. English transl. J. Wiley, Chichester-New York, 1984.

\bibitem{B14} \textsc{P. Baroni:} \emph{Lorentz estimates for obstacle parabolic problems}, Nonlinear Anal., {\bf 96}, (2014), 167-188.

\bibitem{BFM01} \textsc{M. Bildhauer, M. Fuchs, G. Mingione:} \emph{A priori gradient bounds and local $\mathcal{C}^{1, \alpha}$-estimates for (double) obstacle problems under non-standard growth conditions}, Z. Anal. Anwendungen, {\bf 20}, (4), (2001), 959-985.

\bibitem{BDM} \textsc{V. B\"ogelein, F. Duzaar, G. Mingione:} \emph{Degenerate problems with irregular obstacles}, J. Reine Angew. Math., {\bf 650}, (2011), 107-160.

\bibitem{BLS15} \textsc{V. B\"ogelein, T. Lukkari, C. Scheven:} \emph{The obstacle problem for the porous medium equation}, Math. Ann. {\bf 363}, (1) (2015), 455-499.

\bibitem{BK74} \textsc{H. Br\'ezis, D. Kinderlehrer:} \emph{The smoothness of solutions to nonlinear variational inequalities}, Indiana Univ. Math. J., {\bf 23} (1973-1974), 831-844.

\bibitem{BCO} \textsc{S.S. Byun, Y. Cho, J. Ok:} \emph{
Global gradient estimates for nonlinear obstacle problem  with non standard growth,}  Forum Math., {\bf 28}, (4), (2016), 729-747.


\bibitem{C76} \textsc{L. Caffarelli:} \emph{The regularity of elliptic and parabolic free boundaries}, Bull. Amer. Math. Soc., {\bf 82}, (1976), 616-618.

\bibitem{CK80} \textsc{L.A. Caffarelli, D. Kinderlehrer:} \emph{Potential methods in variational inequalities}, J. Anal. Math., {\bf 37} (1980), 285-295.

\bibitem{C84} \textsc{M. Chipot:} \emph{Variational inequalities and flow in porous media}, Springer-Verlag, Berlin and New York, 1984.

\bibitem{C91} \textsc{H. Choe:} \emph{A regularity theory for a general class of quasilinear elliptic partial differential equations and obstacle problems}, Arch. Rational Mech. Anal., {\bf 114}, (4), (1991), 383-394.

\bibitem{CL91} \textsc{H. Choe, J.-L. Lewis:} \emph{On the obstacle problem for quasilinear elliptic equations of $p$ Laplacian type}, SIAM J. Math. Anal., {\bf 22}, (3), (1991), 623-638.


%

\bibitem{CGGP1} \textsc{G. Cupini, F. Giannetti, R. Giova, A. Passarelli di Napoli:} \emph{Higher integrability for minimizers of asymptotically convex integrals with discontinuous coefficients}, Nonlinear Anal., {\bf 154}, (2017), 7-24. 

\bibitem{CGGP2} \textsc{G. Cupini, F. Giannetti, R. Giova, A. Passarelli di Napoli:} \emph{Regularity results for vectorial minimizers of a class
of degenerate convex integrals,} J. Differential Equations {\bf 265}, (2018),  
no. 9, 4375–4416. 

\bibitem{CGGP2} \textsc{G. Cupini, P. Marcellini, E. Mascolo:} \emph{Existence and regularity for elliptic equations under p,q -growth,} Advances in Differential Equations, 19(7/8), (2014), 693-724.

\bibitem{D87} \textsc{F. Duzaar:} \emph{Variational inequalities for harmonic maps,} J. Reine Angew. Math., {\bf 374} (1987), 39-60.

\bibitem{E07} \textsc{M. Eleuteri:} \textit{Regularity
results for a class of obstacle problems}, Appl. Math., {\bf 52}, (2), (2007), 137-169.

\bibitem{EH08} \textsc{M. Eleuteri, J. Habermann:}
\textit{Regularity results for a class of obstacle problems under
non standard growth conditions}, J. Math. Anal. Appl., {\bf 344}, (2), (2008), 1120-1142.

\bibitem{EH10} \textsc{M. Eleuteri, J. Habermann:}
\textit{Calder\'on-Zygmund type estimates for a class of obstacle
problems with $p(x)$ growth}, J. Math. Anal. Appl., {\bf 372}, (1), (2010), 140-161.

\bibitem{EH11} \textsc{M. Eleuteri, J. Habermann:}
\textit{A H\"older continuity result for a class of obstacle
problems under non standard growth conditions}, Math. Nachr., {\bf 284}, No. 11-12, (2011), 1404-1434.

\bibitem{EHL13} \textsc{M. Eleuteri, P. Harjulehto, T. Lukkari:}
\textit{Global regularity and stability of solutions to obstacle problems with nonstandard growth}, Rev. Mat. Complut., {\bf 26}, No. 1, (2013), 147-181.

\bibitem{EMM} \textsc{M. Eleuteri, P. Marcellini, E. Mascolo:} Lipschitz estimates for systems with ellipticity conditions at infinity, \textit{Ann. Mat. Pura Appl.,}  {\bf 195} (5), (2016), 1575-1603.

\bibitem{EMM} \textsc{M. Eleuteri, P. Marcellini, E. Mascolo:} \emph{Regularity for scalar integrals without structure conditions,} Advances in Calculus of Variations, (2018), DOI: https://doi.org/10.1515/acv-2017-0037

\bibitem{EPdN1} \textsc{M. Eleuteri, A. Passarelli di Napoli:}{\ Higher differentiability for solutions to a class of obstacle problems}, Calc. Var. Partial Differential Equations, {\bf 57}: 115 (2018).

\bibitem{EPdN2} \textsc{M. Eleuteri, A. Passarelli di Napoli:}{\ Regularity results for a class of non-differentiable obstacle problems}, Nonlinear Anal,  (2019) \textit{(to appear)}. DOI: 10.1016/j.na.2019.01.024

\bibitem{F64} \textsc{G. Fichera:} \emph{Problemi elastostatici con vincoli unilaterali: il problema di Signorini con ambigue condizioni al contorno,} Atti Accad. Naz. Lincei Mem. Cl. Sci. Fis. Mat. Nat. Sez. Ia, {\bf 7}, (8), (1963-1964), 91-140.

\bibitem{F82} \textsc{A. Friedman:} \emph{Variational principles and free boundary problems}, Wiley, New York, 1982.

\bibitem{F85} \textsc{M. Fuchs:} \emph{Variational inequalities for vector valued functions with non
convex obstacles,} Analysis {\bf 5} (1985), 223-238.

\bibitem{F90} \textsc{M. Fuchs:} \emph{p-harmonic obstacle problems. Part I: Partial regularity theory,} Ann. Mat. Pura Appl. {\bf 156} (1990), 127-158.

\bibitem{F94} \textsc{M. Fuchs:} \emph{Topics in the Calculus of Variations}, Advanced Lectures in Mathematics, Vieweg, (1994).

\bibitem{FG} \textsc{M. Fuchs, L. Gongbao:} \emph{Variational inequalities for energy functionals with nonstandard growth conditions}, Abstr. Appl. Anal., {\bf 3} (1998), 41-64.

\bibitem{FM00} \textsc{M. Fuchs, G. Mingione:} \emph{Full $\mathcal{C}^{1, \alpha}-$regularity for free and constrained local minimizers of elliptic variational integrals with nearly linear growth}, Manuscripta Math., {\bf 102}, (2000), 227-250.

\bibitem{G19} \textsc{C. Gavioli:} \emph{Higher differentiability for a class of obstacle problems with nonstandard growth conditions}, (2019) Forum Mathematicum. 2019, Vol. 31, No. 6, pp. 1501-1516.

\bibitem{G20} \textsc{C. Gavioli:} \emph{A priori estimates for solutions to a class of obstacle problems under p, q- growth conditions.} Journal of Elliptic and Parabolic Equations. 2019, Vol. 5, No. 2, pp. 325-347.

\bibitem{giaquinta}\textsc{M. Giaquinta:}\emph{\ \ Multiple integrals in the
calculus of variations and nonlinear elliptic systems,}{\ \ Annals of
Mathematics Studies, 105. Princeton University Press.}{\ (1983)}.

\bibitem{giusti}\textsc{E. Giusti:}\emph{\ 
Direct methods in the calculus of variations,}{\  World Scientific Publishing Co., Inc., River Edge, NJ}{\ (2003)}.


\bibitem{KS80} \textsc{D. Kinderlehrer, G. Stampacchia:} \emph{An introduction to variational inequalities and their applications}, Academic Press, NY (1980).

\bibitem{L88} \textsc{P. Lindqvist:} \emph{Regularity for the gradient of the solution to a nonlinear obstacle problem with degenerate ellipticity}, Nonlinear Anal., {\bf 12}, (11), (1988), 1245-1255.

\bibitem{LS67} \textsc{J.L. Lions, G. Stampacchia:} \emph{Variational inequalities}, Comm. Pure Appl. Math., {\bf 20}, (1967), 493-519.

\bibitem{MZ86} \textsc{J.H. Michael, W.P. Ziemer:} \emph{Interior regularity for solutions to obstacle problems}, Nonlinear Anal., {\bf 10}, (12), (1986), 1427-1448.


\bibitem{M1}\textsc{P. Marcellini:}{\ Regularity of minimizers of
integrals in the calculus of variations with non standard growth conditions,}\emph{\ Arch. Rational Mech. Anal.}\textbf{\ 105}{\ (1989) 267-284}.


\bibitem{M2}\textsc{P. Marcellini:}{\ Regularity and existence of
solutions of elliptic equations with $p-q$-growth conditions,}\emph{\ J.
Differential Equations}\textbf{\ 90}{\ (1991) 1-30}.

\bibitem{M3}\textsc{P. Marcellini:}{\ Regularity for elliptic equations
with general growth conditions,}\emph{\ J. Differential Equations}\textbf{\
105}{\ (1993) 296-333}.

\bibitem{M4}\textsc{P. Marcellini:}{\ Regularity for some scalar variational problems under general growth conditions,}\emph{\ J. Optim. Theory Appl.}\textbf{\
90, no. 1}{\ (1996) 161-181}.


\bibitem{Moser} \textsc{J. Moser:} \emph{A new proof of De Giorgi's theorem concerning the regularity problem for elliptic differential equations,} Comm. Pure Appl. Math., {\bf 13} (1960), 475-468.

\bibitem{O16-1} \textsc{J. Ok:} \emph{Regularity results for a class of obstacle problems with nonstandard growth}, J. Math. Anal. Appl., {\bf 444}, (2), (2016), 957-979.

\bibitem{O16-2} \textsc{J. Ok:} \emph{Calder\'on-Zygmund estimates for a class of obstacle problems with nonstandard growth}, NoDEA, Nonlinear Diff. Equ. Appl., {\bf 23}, (4), (2016).

\bibitem{O17} \textsc{J. Ok:} \emph{Gradient continuity for nonlinear obstacle problems}, Mediterr. J. Math., {\bf 14}, (1), (2017).


\bibitem{R87} \textsc{J.F. Rodrigues:} \emph{Obstacle problems in mathematical physics}, Elsevier 1987.

\bibitem{S64} \textsc{G. Stampacchia:} \emph{Formes bilineaires coercivitives sur les ensembles convexes,} C.R. Ac. Sci. Paris, {\bf 258}, (1964), 4413-4416.


\end{thebibliography}
\end{document}